\documentclass{amsart}
\usepackage{amssymb}
\usepackage[all]{xy}
\usepackage{slashbox}
\usepackage{rotating}

\begin{document}
\bibliographystyle{plain}

\theoremstyle{definition}
\newtheorem{defn}{Definition}[section]
\theoremstyle{plain}
\newtheorem{thm}[defn]{Theorem}
\newtheorem{prop}[defn]{Proposition}
\newtheorem{cor}[defn]{Corollary}
\newtheorem{lemma}[defn]{Lemma}
\theoremstyle{remark}
\newtheorem{example}[defn]{Example}
\newtheorem{examples}[defn]{Examples}
\newtheorem{rem}[defn]{Remark}

\newcommand{\LCA}{\mathsf{LCAb}}
\newcommand{\FLCA}{\mathsf{FLCAb}}
\newcommand{\TAb}{\mathsf{TAb}}
\newcommand{\Quis}{\mathsf{Quis}}
\newcommand{\Div}{\mathsf{Div}}
\newcommand{\Codiv}{\mathsf{Codiv}}
\newcommand{\I}{\mathsf{I}}
\renewcommand{\P}{\mathsf{P}}
\newcommand{\C}{\mathsf{C}}
\newcommand{\D}{\mathsf{D}}
\newcommand{\K}{\mathsf{K}}
\newcommand{\LH}{\mathcal{LH}}
\newcommand{\Ab}{\mathsf{Ab}}
\newcommand{\dual}{\vee}
\newcommand{\disc}{\mathrm{disc}}
\newcommand{\tors}{\mathrm{tors}}
\newcommand{\toptors}{\mathrm{toptors}}
\newcommand{\naturals}{\mathbb{N}}
\newcommand{\integers}{\mathbb{Z}}
\newcommand{\rationals}{\mathbb{Q}}
\newcommand{\reals}{\mathbb{R}}
\newcommand{\adeles}{\mathbb{A}}
\newcommand{\finadeles}{\mathbb{A}_{\mathrm{fin}}}
\newcommand{\torus}{\mathbb{S}^1}
\newcommand{\field}{\mathbb{F}}
\newcommand{\solenoid}{{\rationals^{\dual}}}
\newcommand{\op}{\mathrm{op}}
\newcommand{\gr}{\mathrm{gr}}
\newcommand{\id}{\mathrm{id}}
\renewcommand{\H}{\mathrm{H}}
\newcommand{\Hom}{\mathrm{Hom}}
\newcommand{\RHom}{\mathrm{RHom}}
\newcommand{\Ext}{\mathrm{Ext}}
\newcommand{\Lotimes}{\otimes^{\mathrm{L}}}
\newcommand{\coker}{\mathrm{coker}}
\newcommand{\im}{\mathrm{im}}
\newcommand{\longto}[1][]{\stackrel{#1}{\longrightarrow}}
\renewcommand{\to}[1][]{\stackrel{#1}{\rightarrow}}
\newcommand{\into}{\hookrightarrow}
\newcommand{\longgets}[1][]{\stackrel{#1}{\longleftarrow}}
\newcommand{\resprod}{\lefteqn{\prod}\!\coprod}
\newcommand{\textresprod}{\lefteqn{\textstyle \prod}\!\coprod}
\newcommand{\adsequence}{[\rationals \to \adeles \to \solenoid]}

\title[Homological algebra with LCA groups]{Homological algebra\\with locally compact abelian groups}
\author[N. Hoffmann]{Norbert Hoffmann}
\address{Norbert Hoffmann\\Mathematisches Institut der Georg-August-Universit\"at\\Bunsenstra{\ss}e 3--5\\
37073 G\"ottingen\\Germany}
\curraddr{School of Mathematics\\Tata Institute of Fundamental Research\\Homi Bhabha Road\\Mumbai 400005\\
India}
\email{hoffmann@uni-math.gwdg.de}

\author[M. Spitzweck]{Markus Spitzweck}
\address{Markus Spitzweck\\Mathematisches Institut der Georg-August-Universit\"at\\Bunsenstra{\ss}e 3--5\\
37073 G\"ottingen\\Germany}
\email{spitz@uni-math.gwdg.de}


\begin{abstract}
  In this article we study \emph{locally compact} abelian (LCA) groups from the
  viewpoint of derived categories, using that
  their category is quasi-abelian in the sense of J.-P. Schneiders. We define a well-behaved
  derived Hom-complex
  with values in the derived category of Hausdorff topological abelian groups.
  Furthermore we introduce a smallness condition for LCA groups and show that such
  groups have a natural tensor product and internal Hom which both admit derived versions.
\end{abstract}
\maketitle

\section*{Introduction}

The aim of this work is to develop further homological algebra techniques for appropriate
\emph{topological} abelian groups, mainly for locally compact abelian (LCA) groups. Already M. Moskowitz
\cite{moskowitz} and also Fulp-Griffith \cite{fulpgriffith} have undertaken to carry over the classical
theory \`{a} la Cartan-Eilenberg; however, their results are limited by the fact that the category $\LCA$ of
LCA groups is not abelian and has neither enough injectives nor enough projectives, as Moskowitz proves.

To solve the first problem, we note that the category $\LCA$ is still \emph{quasi-abelian}, a term introduced
by J.-P. Schneiders \cite{schneiders} which allows to form a derived category with a good deal of the usual
properties. As to the second problem, we substitute injectivity by the weaker notion of divisibility, and
projectivity by the dual notion of codivisibility; it turns out that these are just good enough to derive
the $\Hom$-functor (even its topological version involving the compact-open topology), by first resolving
both variables and then introducing some explicit `correction term' which compensates the fact that the
objects in these resolutions are only `almost acyclic'.

The structure of this paper is as follows: Section \ref{quasiab} introduces the categories of LCA groups and
more generally of abelian Hausdorff groups, shows that they are quasi-abelian and summarizes the consequences
according to \cite{schneiders}: They yield derived categories, and they naturally embed into abelian ones.

Section \ref{types+ranks} contains some structure theory of LCA groups, which may be of independent interest:
We observe that every LCA group has a canonical filtration of length three; this generalises the canonical
torsion subgroup of a discrete abelian group and the dual canonical subgroup of a compact abelian group. We
also introduce some smallness property for LCA groups which we call `finite ranks' and which familiar
examples like $\rationals/\integers$, $\rationals$, $\reals$, $\rationals_p$, $\hat{\integers}$ and $\adeles$
satisfy; their category has a tensor product and internal $\Hom$, as we can prove at the end of the next
section.

Section \ref{divisibility} deals with topological analogues of the fact that divisible abelian groups are
injective in $\Ab$; the corresponding statements in $\LCA$ require much stronger hypotheses, and the proofs
are more involved. We also show topological analogues of the fact that every abelian group can be resolved
by divisible ones.

Finally, section \ref{derived_Hom} contains our construction of a derived $\Hom$-functor; we obtain an exact
bifunctor from the bounded derived category of LCA groups to the bounded derived category of abelian
Hausdorff groups, show that its zeroth cohomology gives the morphism group in the former derived category,
and refine the known fact that $\LCA$ has cohomological dimension $1$. There are also derived versions of
the tensor product and internal $\Hom$ for LCA groups that have finite ranks; we compute the resulting ring
structure on their Grothendieck group $K_0$.

Both authors thank the Mathematical Institute of G\"ottingen University, especially Y. Tschinkel, for their
stimulating interest in this work. The first author also thanks the Tata Institute of Fundamental Research
in Mumbai for its hospitality and support during the preparation of this text.

\section{The quasi-abelian categories $\LCA$ and $\TAb$} \label{quasiab}

This text deals with \emph{locally compact abelian (LCA) groups}, i.\,e.\ abelian topological groups whose
underlying topological space is locally compact, in particular Hausdorff. Standard examples are
the discrete groups $\integers$, $\integers/n$, $\rationals$, $\rationals/\integers$ and $\rationals_p/
\integers_p$; the Lie groups $\reals$ and $\torus := \reals/\integers$; the Pontryagin dual $A^{\dual} :=
\Hom( A, \torus)$ of any LCA group $A$, e.\,g.\ the \emph{solenoid} $\solenoid$; the profinite groups
$\integers_p$ and $\hat{\integers} \cong \prod_p \integers_p$; the additive groups of $p$-adic numbers
$\rationals_p$ and their restricted product, the additive group of finite adeles $\finadeles :=
\textresprod_p (\rationals_p : \integers_p)$; the group of all adeles $\adeles := \reals \oplus \finadeles$.

LCA groups and continuous homomorphisms form an additive category $\LCA$; this is a full additive subcategory
of the additive category $\TAb$ of all abelian Hausdorff groups. If $A$ is an LCA group and $A' \subseteq
A$ is a closed subgroup, then $A'$ and $A/A'$ (endowed with the induced topology) are again LCA groups;
similarly for abelian Hausdorff groups. In particular, each morphism $f:A \to B$ in $\LCA$ or in $\TAb$
has a kernel $\ker( f) := f^{-1}( 0)$ and a cokernel $\coker( f) := B/\overline{f( A)}$. We call $f$ a
\emph{monomorphism} if $\ker( f) = 0$ and an \emph{epimorphism} if $\coker( f)$ = 0. The categories $\LCA$
and $\TAb$ are not abelian. For example, the inclusion morphism $\rationals \to \reals$ is a monomorphism
and an epimorphism, but no isomorphism.
\begin{defn}
  i) A morphism $f: A \to B$ in $\LCA$ or in $\TAb$ is \emph{strict} if the induced monomorphism
  $\bar{f}: A/\ker( f) \to B$ is a closed embedding.

  ii) A complex $A^{\bullet}$ of LCA groups or abelian Hausdorff groups is \emph{strictly exact} if it is
  exact and $\partial: A^n \to A^{n+1}$ is strict for all $n$.
\end{defn}
This definition is taken from Schneiders \cite{schneiders}, who more generally calls a morphism $f$ in an
additive category with kernels and cokernels \emph{strict} if the induced morphism $\bar{f}: \coker( \ker(f))
\to \ker(\coker(f))$ is an isomorphism. In the case of LCA groups, Moskowitz \cite{moskowitz} and Armacost
\cite{armacost} call such morphisms \emph{proper}.

If a morphism $f: A \to B$ in $\LCA$ is strict, then $f^*: B^{\dual} \to A^{\dual}$ also is. A
monomorphism (resp. epimorphism) $f:A \to B$ in $\LCA$ or in $\TAb$ is strict if and only if $f$ is closed
(resp. open); cf. \cite[Thms. 5.26 and 5.27]{hewittross} for the latter. In particular, the composition of
two strict monomorphisms (resp. epimorphisms) is again strict; cf. also \cite[Prop. 1.1.7]{schneiders}.
By contrast, the composition of two arbitrary strict morphisms is in general not strict.

\begin{prop} \label{is_quasiabelian}
  If the commutative diagram in $\LCA$ or in $\TAb$
  \begin{equation*} \xymatrix{
    A \ar[r]^{\alpha} \ar[d]_f & A' \ar[d]^{f'}\\
    B \ar[r]_{\beta}           & B'
  } \end{equation*}

  i) is a pushout square and $f$ is a strict monomorphism, then $f'$ also is.

  ii) is a pullback square and $f'$ is a strict epimorphism, then $f$ also is.
\end{prop}
\begin{proof}
  i) If $f$ is closed and injective, then
  \begin{equation*}
    F := \begin{pmatrix} 1 & -\alpha\\0 & f \end{pmatrix}
       = \begin{pmatrix} 1 &       0\\0 & f \end{pmatrix}
         \begin{pmatrix} 1 & -\alpha\\0 & 1 \end{pmatrix}: A' \oplus A \longto A' \oplus B
  \end{equation*}
  also is. But $f'$ is obtained from $F$ by dividing out the closed subgroup $\{0\} \times A$ and its image
  under $F$, so $f'$ is closed and injective as well.

  ii) If $f'$ is open and surjective, then
  \begin{equation*}
    F' := \begin{pmatrix} f' & -\beta\\0 & 1 \end{pmatrix}
        = \begin{pmatrix} 1  & -\beta\\0 & 1 \end{pmatrix}
          \begin{pmatrix} f' &      0\\0 & 1 \end{pmatrix}: A' \oplus B \longto B' \oplus B
  \end{equation*}
  also is. But $f$ is just the restriction of $F'$ to the inverse image of $\{0\} \times B$, hence open and
  surjective as well.
\end{proof}
This proposition means that the categories $\LCA$ and $\TAb$ are \emph{quasi-abelian} in the sense of J.-P.
Schneiders \cite{schneiders}, so his results apply as follows:

Starting from the categories of bounded complexes $C^b( \LCA) \subseteq \C^b( \TAb)$ and identifying chain
homotopic morphisms, we obtain as usual triangulated categories $\K^b( \LCA) \subseteq \K^b( \TAb)$. A
morphism in $\K^b( \LCA)$ or in $\K^b( \TAb)$ is called a \emph{strict quasi-isomorphism} if its mapping
cone is strictly exact. \cite[Prop. 1.2.14]{schneiders} states that the class of strictly exact complexes in
$\K^b( \LCA)$ or in $\K^b( \TAb)$ is \emph{stable under extensions} in the sense of \cite[IV.2.10]{manin}.
This implies that the class $\Quis$ of strict quasi-isomorphisms is \emph{localizing}
\cite[Def. III.2.6]{manin} and also \emph{compatible with the triangulation} \cite[IV.2.1]{manin}; thus one
obtains \emph{derived categories} $\D^b(\LCA) := \K^b( \LCA)/\Quis$ and $\D^b( \TAb) := \K^b( \TAb)/\Quis
$ which are again triangulated due to \cite[Thm. IV.2.2]{manin}. See also \cite[Thm. 2.1.8]{neeman} for this
quotient construction. We remark that a priori it is not clear that these quotient constructions yield
genuine categories in the sense that the Hom object between two fixed objects is (isomorphic to) a set. This
will follow in the case of $\D^b(\LCA)$ from proposition \ref{RHom_lim}.i, in the case of $\D^b(\TAb)$ we do
not know it and just note that there do not arise set-theoretic difficulties in
allowing proper classes as Hom objects.

\cite[Def. 1.2.18]{schneiders} endows $\D^b(\LCA)$ and $\D^b(\TAb)$ with (`left') t-structures, yielding
abelian categories $\LH( \LCA)$ and $\LH( \TAb)$ as their heart. These are `abelian envelopes' in the sense
that the natural functors $I: \LCA \hookrightarrow \LH( \LCA)$ and $I: \TAb \hookrightarrow \LH(\TAb)$
are fully faithful \cite[Cor. 1.2.28]{schneiders}; their essential image is stable under extensions and
subobjects, and it contains the cokernel of $I(f)$ for a morphism $f$ in $\LCA$ or $\TAb$ if and only if
$f$ is strict \cite[Prop. 1.2.29]{schneiders}. According to \cite[Prop. 1.2.32]{schneiders}, $I$ induces an
equivalence of derived categories.

We will sometimes identify groups with discrete topological groups; this defines a fully faithful embedding
$\Ab \hookrightarrow \TAb\hookrightarrow \LH(\TAb)$. This embedding has an exact left inverse $\LH(
\TAb) \to \Ab$ which we denote by $A \mapsto A_{\disc}$ and which is given as follows: It sends each object
$A$ of $\TAb$ to its underlying discrete group $A_{\disc}$; since this preserves kernels of arbitrary and
cokernels of strict morphisms, it does induce an exact functor $\LH( \TAb) \to \Ab$ by
\cite[Prop. 1.2.34]{schneiders}.

\section{Types of LCA groups and finite ranks} \label{types+ranks}

Recall that an LCA group $A$ is a \emph{topological $p$-group} (resp. a \emph{topological torsion group}) if
$\lim_{n \to \infty} p^n a = 0$ (resp. $\lim_{n \to \infty} n! a = 0$) holds for all $a \in A$. According to
the Braconnier-Vilenkin theorem \cite[Thm. 3.13]{armacost}, $A$ is a topological torsion group if and only if
$A \cong \textresprod_p (A_p : U_p)$ is a restricted product of topological $p$-groups $A_p$ with respect to
compact open subgroups $U_p \subseteq A_p$. $A_p$ is uniquely determined as a closed subgroup of $A$; it is
called the \emph{$p$-component} of $A$.
\begin{defn}
  An LCA group $A$ is

  i) of \emph{type $\integers$} if $A$ is discrete and torsionfree,

  ii) of \emph{type $\torus$} if $A$ is compact and connected,

  iii) of \emph{type $\adeles$} if $A \cong A_{\reals} \oplus A_{\toptors}$ with $A_{\reals} \cong \reals^n$
  for some $n$ and $A_{\toptors}$ a topological torsion group.

\end{defn}
Note that the direct sum decomposition $A \cong A_{\reals} \oplus A_{\toptors}$ in iii is unique and
functorial, because $\Hom( A_{\reals}, A_{\toptors}) = 0 = \Hom( A_{\toptors}, A_{\reals})$.

$A$ is a topological $p$-group (resp. a topological torsion group, resp. of type $\adeles$) if and only if
its Pontryagin dual $A^{\dual}$ is \cite[Cor. 2. 13 and Cor. 3.7]{armacost}; $A$ is of type $\torus$ if and
only if $A^{\dual}$ is of type $\integers$ \cite[Cor. 4 to Thm. 30]{morris}.

\begin{prop} \label{types}
  i) Every LCA group $A$ has a unique chain of closed subgroups $0 \subseteq A_{\torus} \subseteq
  F_{\integers} A \subseteq A$ such that $A_{\torus}$, $A_{\adeles} := F_{\integers} A/A_{\torus}$ and
  $A_{\integers} := A/F_{\integers} A$ are of type $\torus$, $\adeles$ and $\integers$, respectively.

  ii) $f( A_{\torus}) \subseteq B_{\torus}$ and $f( F_{\integers} A) \subseteq F_{\integers} B$ for all
  morphisms $f: A \to B$ in $\LCA$.
\end{prop}
\begin{proof}
  i) Existence: By the structure theorem for LCA groups \cite[Thm. 24.30]{hewittross}, $A \cong \reals^n
  \oplus A'$ where $A'$ has a compact open subgroup. Without loss of generality, we may thus assume that $A$
  itself has a compact open subgroup $U \subseteq A$.

  Let $A_{\torus} \subseteq U$ be the connected component of $0$; this closed subgroup of $A$ has type
  $\torus$. Let $F_{\integers} A \subseteq A$ be the inverse image of the torsion in the discrete group
  $A/U$; this is an open subgroup such that $A/F_{\integers} A$ has type $\integers$.

  The open subgroup $U/A_{\torus}$ of $A_{\adeles} := F_{\integers}A/A_{\torus}$ is compact
  and totally disconnected, hence profinite; the quotient $F_{\integers}A/U$ is a discrete
  torsion group. Consequently, $A_{\adeles}$ is topological torsion, in particular of type $\adeles$.

  Uniqueness and ii follow from the observation that all morphisms $A_{\torus} \to A_{\adeles}$, $A_{\torus}
  \to A_{\integers}$ and $A_{\adeles} \to A_{\integers}$ vanish if $A_{\torus}$, $A_{\adeles}$ and
  $A_{\integers}$ are of type $\torus$, $\adeles$ and $\integers$.
\end{proof}

For a general LCA group $A$, we put $A_{\reals} := (A_{\adeles})_{\reals}$ and $A_{\toptors} := (A_{\adeles}
)_{\toptors}$; we write $A_p$ for the $p$-components of the latter. For a morphism $f: A \to B$ in $\LCA$, we
denote by $f_?:A_? \to B_?$ the induced morphisms for $? \in \{\torus,\adeles,\reals,\toptors,p,\integers\}$.

\begin{lemma} \label{types_extension}
  If $0 \to A \to B \to[\pi] C \to 0$ is a strictly exact sequence of LCA groups in which $A$ and $C$ are
  both of type $\torus$ (resp. both of type $\adeles$, resp. both of type $\integers$), then $B$ is also of
  type $\torus$ (resp. $\adeles$, resp. $\integers$).
\end{lemma}
\begin{proof}
  If $A$ and $C$ are discrete (resp. torsionfree), then $B$ also is. This proves the case of
  type $\integers$. The case of type $\torus$ follows by duality.

  If $A$ and $C$ are of type $\adeles$, then $\pi( B_{\torus}) \subseteq C_{\torus} = 0$ by proposition
  \ref{types}.ii, so $B_{\torus} \subseteq A$ and hence $B_{\torus} \subseteq A_{\torus} = 0$ by proposition
  \ref{types}.ii again. This shows $B_{\torus} = 0$ and by duality also $B_{\integers} = 0$; hence $B$
  is of type $\adeles$.
\end{proof}

Given a strictly exact sequence of LCA groups, what can we say about their types? The following sequences and
their duals are examples of mixed types:
\begin{equation*} \begin{array}{ccccccccc}
  0 & \longto & \integers & \longto[\cdot n] & \integers & \longto & \integers/n
    & \longto & 0\\[1ex]
  0 & \longto & \integers & \longto & \rationals & \longto & \rationals/\integers
    & \longto & 0\\
  0 & \longto & \integers & \longto[\binom{p^n}{1}] & \integers \oplus \integers_p
    & \longto[(-1 \; p^n)] & \integers_p & \longto & 0\\[1ex]
  0 & \longto & \integers & \longto & \reals & \longto & \torus & \longto & 0\\
  0 & \longto & \integers & \longto[\binom{n}{-1}] & \reals \oplus \integers/n
    & \longto & \torus & \longto & 0\\[1ex]
  0 & \longto & \rationals & \longto & \reals \oplus \rationals/\integers
    & \longto & \torus & \longto & 0\\[1ex]
  0 & \longto & \rationals & \longto & \adeles & \longto & \solenoid & \longto & 0
\end{array} \end{equation*}

\begin{prop} \label{types_sequence}
  Given a strictly exact sequence of LCA groups
  \begin{equation*}
    0 \longto A \longto[\iota] B \longto[\pi] C \longto 0,
  \end{equation*}
  there are unique closed subgroups $0 = F_0 B \subseteq F_1 B \subseteq \ldots \subseteq F_6 B \subseteq F_7
  B = B$ such that the induced subgroups $F_n A := \iota^{-1}( F_n B)$ and $F_n C := \pi( F_n B)$ are also
  closed and the induced sequences
  \begin{equation*}
    0 \longto        \gr_n A := \frac{F_n A}{F_{n-1} A}
      \longto[\iota] \gr_n B := \frac{F_n B}{F_{n-1} B}
      \longto[\pi]   \gr_n C := \frac{F_n C}{F_{n-1} C} \longto 0
  \end{equation*}
  are strictly exact of the following types:
  \begin{equation*} \begin{array}{ccccccccc@{\qquad}c}
    0 & \longto & \text{type } \torus    & \longto & \text{type } \torus   
      & \longto & 0                      & \longto & 0 & n = 1\\
    0 & \longto & \text{profinite}       & \longto & \text{type } \torus
      & \longto & \text{type } \torus    & \longto & 0 & n = 2\\
    0 & \longto & \text{type } \adeles   & \longto & \text{type } \adeles
      & \longto & 0                      & \longto & 0 & n = 3\\
    0 & \longto & \text{type } \integers & \longto & \text{type } \adeles
      & \longto & \text{type } \torus    & \longto & 0 & n = 4\\
    0 & \longto & 0                      & \longto & \text{type } \adeles
      & \longto & \text{type } \adeles   & \longto & 0 & n = 5\\
    0 & \longto & \text{type } \integers & \longto & \text{type } \integers
      & \longto & \text{torsion}         & \longto & 0 & n = 6\\
    0 & \longto & 0                      & \longto & \text{type } \integers
      & \longto & \text{type } \integers & \longto & 0 & n = 7
  \end{array} \end{equation*}
\end{prop}
\begin{proof}
  Uniqueness: Due to lemma \ref{types_extension}, $F_2 B$ has to be of type $\torus$, $F_5 B \big/ F_2 B$ of
  type $\adeles$, and $B \big/ F_5 B$ of type $\integers$. Hence $F_2 B = B_{\torus}$ and $F_5 B = F_{
  \integers} B$ are uniquely determined by proposition \ref{types}.i. Similarly, $F_1 A = A_{\torus}$ and
  $F_3 A = F_{\integers} A$; thus $F_1 B = \iota( A_{\torus})$ and $F_3 B = F_2 B + \iota( F_{\integers} A)$
  since we require $\gr_1 C = \gr_3 C = 0$. By the dual argument, $F_4 B$ and $F_6 B$ are also uniquely
  determined.

  Existence: We put $F_2 B := B_{\torus}$ and $F_2 A := \iota^{-1}(F_2 B)$, $F_2 C := \pi(F_2 B)$. These are
  compact, in particular closed, subgroups, and $0 \to F_2 A \to F_2 B \to F_2 C \to 0$ is an exact sequence
  by construction, so it is even strictly exact. Using the $3 \times 3$-lemma in the abelian envelope
  $\LH( \LCA)$ of $\LCA$, it follows that the induced sequence $0 \to A/F_2 A \to B/F_2 B \to C/F_2 C \to 0$
  is also strictly exact. Thus it suffices to prove the proposition for both sequences separately, i.\,e.\ we
  may assume without loss of generality $B_{\torus} = B$ or $B_{\torus} = 0$. Applying the dual argument in
  the latter case, we we may assume that $B$ is of type $\torus$ or of type $\adeles$ or of type $\integers$.

  If $B$ is of type $\integers$, then so is $A$, and $C$ is discrete. In this case, $F_5 B := 0$ and $F_6 B
  := \pi^{-1} (F_{\integers} C)$, the inverse image of the torsion in $C$, defines a chain of subgroups in
  $B$ with the required properties. This proves the case of type $\integers$; the case of type $\torus$
  follows by duality.

  Suppose that $B$ is of type $\adeles$. Then $A_{\torus} \subseteq B_{\torus} = 0$ by proposition
  \ref{types}.ii, so $A_{\adeles}$ is a closed subgroup of $A$. The restriction $\iota: A_{\adeles} \to B$
  automatically decomposes into morphisms $\iota_{\reals}: A_{\reals} \to B_{\reals}$ and $\iota_{\toptors}:
  A_{\toptors} \to B_{\toptors}$ which are closed embeddings because $\iota$ is. Thus $B/\iota( A_{\adeles})
  \cong \coker( \iota_{\reals}) \oplus \coker( \iota_{\toptors})$ is of type $\adeles$. So it suffices to
  prove the claim for the sequence $0 \to A/A_{\adeles} \to B/\iota( A_{\adeles}) \to C \to 0$; in other
  words, we may additionally assume that $A$ is of type $\integers$. Using the dual argument, we may also
  assume that $C$ is of type $\torus$. In this situation, simply $F_3 B = 0$ and $F_4 B = B$ does the trick.
\end{proof}

\begin{defn}
  An LCA group $A$ has

  i) \emph{finite $\integers$-rank} if the real vector space $\Hom( A, \reals)$ has finite dimension,

  ii) \emph{finite $\torus$-rank} if the real vector space $\Hom( \reals, A)$ has finite dimension,

  iii) \emph{finite $p$-rank} if $p \cdot \_: A \to A$ is strict with finite kernel and cokernel.
\end{defn}

\begin{defn}
  An LCA group $A$ has \emph{finite ranks} if $A$ has finite $\integers$-rank, finite $\torus$-rank and
  finite $p$-rank for all prime numbers $p$. $\FLCA \subseteq \LCA$ denotes the full additive subcategory
  consisting of all LCA groups that have finite ranks.
\end{defn}

\begin{prop} \label{pure_ranks}
  An LCA group $A$ has finite ranks if

  i) $A$ is of type $\integers$ and has finite $\integers$-rank or

  ii) $A$ is of type $\torus$ and has finite $\torus$-rank or

  iii) $A$ is of type $\adeles$ and has finite $p$-rank for all $p$.
\end{prop}
\begin{proof}
  i) Any LCA group $A$ of type $\integers$ has finite $\torus$-rank because $\Hom(\reals,A) = 0$. If $A$ also
  has finite $\integers$-rank, then $\dim_{\rationals} (A \otimes_{\integers} \rationals) =: d < \infty$; any
  given $d+1$ elements $a_1, \ldots, a_{d+1} \in A$ are thus contained in a subgroup $\integers^d \cong A'
  \subseteq A$ and hence linearly dependent in $A'/pA'$, a fortiori in $A/pA$. Thus $\dim_{\field_p}(A/pA)
  \leq d < \infty$; since $p \cdot \_: A \to A$ is automatically strict and injective, $A$ has finite
  $p$-rank. This proves i; ii follows by duality.

  iii) Any $A$ of type $\adeles$ has finite $\torus$-rank since $\Hom( \reals, A) \cong \Hom(\reals, A_{
  \reals}) \cong A_{\reals}$ is finite-dimensional. The dual argument shows that $A$ also has finite
  $\integers$-rank.
\end{proof}

\begin{lemma} \label{p_structure}
  A topological $p$-group $A \in \LCA$ has finite $p$-rank

  i) if $A$ is discrete and ${}_pA := \{a \in A: pa = 0\}$ is finite.

  ii) if $A$ is compact and $A/pA$ is finite.

  iii) if and only if $A \cong A_1 \oplus \ldots \oplus A_r$ with $A_i \cong \integers/p^{n_i}$ or
  $\rationals_p/\integers_p$ or $\integers_p$ or $\rationals_p$.
\end{lemma}
\begin{proof}
  $\integers/p^{n_i}$, $\rationals_p/\integers_p$, $\integers_p$ and $\rationals_p$ have finite $p$-rank;
  the `if' part of iii follows.

  If $A$ is discrete and ${}_pA$ is finite, then $A \cong (\rationals_p/\integers_p)^n \oplus A'$ for some
  finite group $A'$ by \cite[Ch. III, Thm. 19.2 and Ex. 19]{fuchs}; this implies i and by duality also ii.

  Let $A$ have finite $p$-rank, and let $A_{\tors} \subseteq A$ be its torsion subgroup. Using
  \cite[Thm. 2.12]{armacost}, we can find an open subgroup $U \subseteq A$ with ${}_pA \cap U = \{0\}$; 
  this implies $A_{\tors} \cap U = \{0\}$, proving that $A_{\tors}$ is closed in $A$ and discrete. Due to
  \cite[Prop. 6.21]{armacost}, every direct summand $\rationals_p/\integers_p$ of $A_{\tors}$ is even a
  topological direct summand of $A$; splitting them off, we may assume that $A_{\tors}$ is finite. Applying
  the same to $B := A^{\dual}$, we may additionally assume that $B_{\tors}$ is finite.

  Then $A' := (B/B_{\tors})^{\dual} \subseteq A$ is an open subgroup of finite index such that $p \cdot \_:
  A' \to A'$ is strict and surjective. In particular, $A' \cap A_{\tors} = 0$ since the latter is finite;
  hence $A'$ is a locally compact topological vector space over $\rationals_p$, so $A' \cong \rationals_p^n$
  by \cite[Ch. I, \S 2]{bourbaki}. Dually, $B' := (A/A_{\tors})^{\dual} \subseteq B$ satisfies $B' \cap B_{
  \tors} = 0$; this means $A_{\tors} + A' = A$. Altogether, we obtain $A = A' \oplus A_{\tors}$, proving iii.
\end{proof}

\begin{rem}
  Part iii of the previous lemma implies that $\rationals_p$ is both injective and projective among
  topological $p$-groups that have finite ranks: It is easy to see that a morphism $A_1 \oplus \ldots \oplus
  A_r \to \rationals_p$ with $A_i \cong \integers/p^{n_i}$ or $\rationals_p/\integers_p$ or $\integers_p$ or
  $\rationals_p$ can only be surjective if its restriction to one summand is an isomorphism.

  By contrast, $\rationals_p$ is neither injective nor projective among all topological $p$-groups; an
  example of a non-split strict epimorphism of topological $p$-groups onto $\rationals_p$ is $\textresprod_{n
  \in \naturals} (\integers_p: p^{2n} \integers_p) \twoheadrightarrow \rationals_p$, $(a_n) \mapsto \sum_n
  p^{-n} a_n$.
\end{rem}

\begin{prop} \label{fin_rank_ext}
  In a strictly exact sequence $0 \to A \to B \to C \to 0$ of LCA groups, $B$ has finite ranks if and only if
  both $A$ and $C$ have finite ranks.
\end{prop}
\begin{proof}
  1) The sequence $0 \to \Hom( \reals, A) \to \Hom( \reals, B) \to \Hom( \reals, C) \to 0$ is exact by
  \cite[Thm. 3.2]{moskowitz}, so $B$ has finite $\torus$-rank if and only if both $A$ and $C$ have finite
  $\torus$-rank. The same holds for their $\integers$-rank by duality.

  2) Note that $A$ has finite $p$-rank if and only if the kernel ${}_pA$ and the cokernel $A/pA$ of $p \cdot
  \_: A \to A$ in $\LH( \LCA)$ are finite groups, i.\,e.\ are in the essential image of the finite groups
  under the embedding $\LCA \hookrightarrow \LH( \LCA)$; similarly for $B$ and $C$. This essential image is
  stable under extensions and subobjects according to \cite[Prop. 1.2.29]{schneiders}; it is also stable
  under quotients because every monomorphism in $\LCA$ into a discrete group is strict. The snake lemma in
  the abelian category $\LH( \LCA)$ yields an exact sequence
  \begin{equation} \label{HomExt}
    0 \longto {}_pA \longto {}_pB \longto {}_pC \longto[\delta] A/pA \longto B/pB \longto C/pC \longto 0
  \end{equation}
  from which we see that $B$ has finite $p$-rank if both $A$ and $C$ have; together with step 1, this proves
  the `if' part of the proposition.

  3) Suppose that $B$ is discrete and has finite ranks. Then $A$ is also discrete; its torsion subgroup
  $A_{\tors}$ has finite $p$-rank by proposition \ref{p_structure}.i because ${}_pA \subseteq {}_pB$ is
  finite, and $A/A_{\tors}$ has finite $p$-rank by proposition \ref{pure_ranks}.i because it has finite
  $\integers$-rank by step 1. Hence $A$ has finite $p$-rank according to step 2; now the exact sequence
  \eqref{HomExt} implies that $C$ also has finite $p$-rank. This proves the `only if' part for discrete $B$;
  it follows by duality for compact $B$. If $B \cong \reals^n$, then $A \cong \reals^a \oplus \integers^b$
  and $C \cong \reals^{n-a-b} \oplus (\torus)^b$ by \cite[Thm. 9.11]{hewittross}, so both have finite ranks
  as well.

  4) Now let $B$ be any LCA group having finite ranks. By the structure theory of LCA groups, we can find
  closed subgroups $B' \subseteq B'' \subseteq B$ such that $B'$ is compact, $B''/B' \cong B_{\reals}$ and
  $B''$ is open. The sequence $0 \to B' \to B \to B/B' \to 0$ yields an exact sequence like \eqref{HomExt}
  with connecting morphism $\delta: {}_p(B/B') \to B'/pB'$ in $\LH( \LCA)$; here ${}_p(B/B')$ is a discrete
  group because $B/B''$ is, $B'/pB'$ is a compact group because $B'$ is, and $\delta$ is strict with finite
  kernel and cokernel because $B$ has finite $p$-rank. Hence ${}_p(B/B')$ and $B'/pB'$ are finite; as ${}_pB'
  \subseteq {}_pB$ is also finite, $B'$ has finite ranks. By duality, $B/B''$ also has.
  
  Let $A' \subseteq A'' \subseteq A$ be the inverse images of $B' \subseteq B'' \subseteq B$, and let $C'
  \subseteq C'' \subseteq C$ be their images. These are closed subgroups because $A'$, $B'$, $C'$ are compact
  and $A''$, $B''$, $C''$ are open; the induced exact sequences $0 \to A' \to B' \to C' \to 0$ and $0 \to
  A/A'' \to B/B'' \to C/C'' \to 0$ are strict for the same reason. By the $3 \times 3$-lemma in the abelian
  category $\LH( \LCA)$, $0 \to A''/A' \to B''/B' \to C''/C' \to 0$ is also strictly exact. Now step 3 shows
  that $A'$, $A''/A'$, $A/A''$ and $C'$, $C''/C'$, $C/C''$ have finite ranks; hence $A$ and $C$ also have
  finite ranks by step 2.
\end{proof}

\begin{cor}
  The category $\FLCA$ is quasi-abelian, and the inclusion functors $\FLCA \hookrightarrow \LCA
  \hookrightarrow \TAb$ preserve kernels and cokernels.
\end{cor}

\begin{rem} \label{finiteness}
  i) If $A^{\bullet} \in \C^b( \LCA)$ is strictly exact, then the discrete group $\H^n( A_{\integers}^{
  \bullet})$ has finite $\integers$-rank for all $n$.

  ii) If $0 \to A \to B \to C \to 0$ is a strictly exact sequence of LCA groups with $A$ of type $\integers$,
  $B$ of type $\adeles$ and $C$ of type $\torus$, then they all have finite ranks.
\end{rem}
\begin{proof}
  i) The strictly exact sequence of complexes $0 \to F_{\integers} A^{\bullet} \to A^{\bullet} \to A_{
  \integers}^{\bullet} \to 0$ yields a long exact sequence in $\LH( \LCA)$ which implies $\H^n( A_{\integers
  }^{\bullet}) \cong \H^{n+1}( F_{\integers} A^{\bullet})$. In particular, the latter is in the essential
  image of $\LCA \hookrightarrow \LH( \LCA)$, so there are closed subgroups $B \subseteq Z \subseteq F_{
  \integers} A^{n+1}$ with $\H^n( A_{\integers}^{\bullet}) \cong Z/B$. Since $\reals$ is injective and
  projective in $\LCA$ by \cite[Thms. 3.2 and 3.3]{moskowitz},
  \begin{equation*}
    \dim \Hom( Z/B, \reals) \leq \dim \Hom( Z, \reals) \leq \dim \Hom( F_{\integers}
    A^{n+1}, \reals) = \dim A^{n+1}_{\reals} < \infty.
  \end{equation*}

  ii) $A$ has finite ranks by i and proposition \ref{pure_ranks}.i; dually, $C$ has finite ranks as well.
  Hence $B$ also has according to proposition \ref{fin_rank_ext}.
\end{proof}

\section{Topological $\Hom$ and divisibility} \label{divisibility}

We usually endow the group $\Hom( A, B)$ of continuous homomorphisms between LCA groups $A$ and $B$ with the
compact-open topology; this turns it into an abelian Hausdorff group. The canonical bijection $\Hom( A, B)
\to \Hom( B^{\dual}, A^{\dual})$, $f \mapsto f^*$ is a topological isomorphism according to
\cite[Cor. 2 to Thm. 4.2]{moskowitz}.

The bifunctor $\Hom: \LCA^{\op} \times \LCA \to \TAb$ is clearly additive in both variables. It is easy to
check directly that $\Hom( A, \_): \LCA \to \TAb$ is left exact, i.\,e.\ preserves kernels; by duality,
$\Hom( \_, B): \LCA^{\op} \to \TAb$ is also left exact, i.\,e.\ transforms cokernels in $\LCA$ to kernels
in $\TAb$.

Given $A^{\bullet}, B^{\bullet} \in \C^b( \LCA)$, we denote by $\Hom^{\bullet}( A^{\bullet}, B^{\bullet})
\in \C^b( \TAb)$ the usual total complex of the $\Hom$-double complex; a special case is the
\emph{dual complex} $(A^{\bullet})^{\dual} := \Hom^{\bullet}( A^{\bullet}, \torus) \in \C^b( \LCA)$.

\begin{defn}
  i) An LCA group $A$ is \emph{divisible} (resp. \emph{strictly divisible}) if the map $n \cdot \_: A \to A$
  is surjective (resp. strict and surjective) for all $n \in \naturals$.

  ii) An LCA group $A$ is \emph{codivisible} if $A^{\dual}$ is divisible.
\end{defn}

\begin{defn}
  $\Div \subseteq \LCA \supseteq \Codiv$ and $\I \subseteq \FLCA \supseteq \P$ denote the full additive
  subcategories given by the following object classes:
  \begin{align*}
    \Div   &= \{ D \in \LCA  \big| D \text{ divisible}\}, &
    \I     &= \{ I \in \FLCA \cap \Div \big| I_{\integers} = 0\},\\
    \Codiv &= \{ C \in \LCA  \big| C \text{ codivisible}\}, &
    \P     &= \{ P \in \FLCA \cap \Codiv \big| P_{\torus} = 0\}.
  \end{align*}
\end{defn}
Note that $\Div \subseteq \LCA$ and $\I \subseteq \FLCA$ are stable under taking quotients; dually, $\Codiv
\subseteq \LCA$ and $\P \subseteq \FLCA$ are stable under taking closed subgroups. This implies in particular
that codivisible LCA groups are torsionfree.

\begin{prop} \label{div_disc}
  If $\iota: A' \hookrightarrow A$ is an open embedding of LCA groups and $D$ is a divisible (resp. strictly
  divisible) LCA group, then $\iota^*: \Hom(A,D) \to \Hom(A',D)$ is surjective (resp. strict and surjective).
\end{prop}
\begin{proof}
  Given a morphism $f': A' \to D$ in $\LCA$, it can be extended to a morphism $f: A \to D$ in $\Ab$ because
  $D$ is an injective object in $\Ab$. But $f$ is automatically continuous: Its restriction to every coset
  modulo $A'$ is continuous because $f'$ is, and these cosets form an open covering of $A$. This shows that
  $\iota^*$ is surjective.

  If $A/A' \cong \integers$, then $\iota$ has a left inverse, so $\iota^*$ is strict. If $A/A' \cong
  \integers/n$ for some $n$, then we can choose an element $a \in A$ whose image generates $A/A'$; this
  defines a pushout square
  \begin{equation*} \xymatrix{
    \integers \ar[r]^{\cdot n} \ar[d] & \integers \ar[d]^{\cdot a}\\ A' \ar[r]^{\iota} & A
  } \end{equation*}
  in $\LCA$. The left exact functor $\Hom( \_, D)$ maps this to a pullback square in $\TAb$; if $D$ is
  strictly divisible, then $\iota^*$ is thus strict and surjective.

  Since the composition of strict epimorphisms is again strict, the proposition is now proved whenever $A/A'$
  is finitely generated. In general, the definition of the compact-open topology yields topological
  isomorphisms
  \begin{equation*}
    \Hom( A,    D) \cong \varprojlim \Hom( \tilde{A},    D) \quad\text{and}\quad
    \Hom( A/A', D) \cong \varprojlim \Hom( \tilde{A}/A', D)
  \end{equation*}
  where both limits are taken over all open subgroups $A' \subseteq \tilde{A} \subseteq A$ such that
  $\tilde{A}/A'$ is finitely generated; now the following lemma completes the proof.
\end{proof}

\begin{lemma} \label{limits}
  Let a sequence of filtered projective systems in $\TAb$
  \begin{equation*}
    0 \longto (A_j)_{j \in J} \longto (B_j)_{j \in J} \longto (C_j)_{j \in J} \longto 0
  \end{equation*}
  be strictly exact on each level $j \in J$. Then the induced sequence of limits
  \begin{equation*}
    0 \longto A := \varprojlim A_j \longto B := \varprojlim B_j \longto C := \varprojlim C_j \longto 0
  \end{equation*}
  is also strictly exact in $\TAb$ if the projections $b_j: B \to B_j$, $c_j: C \to C_j$ and the natural
  map $\ker( b_j) \to \ker( c_j)$ are surjective for all $j$.
\end{lemma}
\begin{proof}
  $A$ is clearly the kernel of the induced map $\pi: B \to C$, so we just have to check that $\pi$ is
  surjective and open. Surjectivity follows from the 5-lemma in $\Ab$, applied to the diagram
  \begin{equation*} \xymatrix{
    0 \ar[r] & \ker( b_j) \ar@{->>}[d] \ar[r] & B \ar[d]^{\pi} \ar[r]^{b_j} & B_j \ar[d]^{\pi_j} \ar[r] & 0\\
    0 \ar[r] & \ker( c_j)              \ar[r] & C              \ar[r]_{c_j} & C_j                \ar[r] & 0;
  } \end{equation*}
  from this diagram also follows that $\pi( b_j^{-1}( U_j))$ equals $c_j^{-1}( \pi_j( U_j))$ and is thus a
  neighborhood of zero for every neighborhood of zero $U_j \subseteq B_j$. But the $b_j^{-1}( U_j)$ form a
  neighborhood base of zero in $B$; this shows that $\pi$ is indeed open.
\end{proof}

\begin{cor} \label{div_split}
  i) For $D \in \Div$, the sequence $0 \to F_{\integers} D \to D \to D_{\integers} \to 0$ splits.

  ii) For $C \in \Codiv$, the sequence $0 \to C_{\torus} \to C \to C/C_{\torus} \to 0$ splits.
\end{cor}
\begin{proof}
  $F_{\integers}D$ is divisible because $D$ is and $D_{\integers}$ is torsionfree. Proposition \ref{div_disc}
  allows to extend the identity on $F_{\integers} D$ to a morphism $D \to F_{\integers} D$ that splits the
  sequence in i. The sequence in ii also splits by duality.
\end{proof}

\begin{prop} \label{acyclic}
  Suppose that $C^{\bullet} = [0 \to C^1 \to C^2 \to C^3 \to 0] \in \C^b( \Codiv)$ and $D^{\bullet} = [0 \to
  D^1 \to D^2 \to D^3 \to 0] \in \C^b( \Div)$ are short strictly exact complexes; let furthermore
  $C \in \Codiv$ and $D \in \Div$ be given.

  i) If $D^n_{\integers} = 0$ for all $n$, then $\Hom^{\bullet}( C, D^{\bullet}) \in \C^b( \TAb)$ is
  strictly exact.

  ii) If $C_{\torus} = 0$, then $\Hom^{\bullet}( C, \adsequence) \in \C^b( \TAb)$ is strictly exact.

  iii) If $C^n_{\torus} = 0$ for all $n$, then $\Hom^{\bullet}( C^{\bullet}, D) \in \C^b( \TAb)$ is
  strictly exact.

  iv) If $D_{\integers} = 0$, then $\Hom^{\bullet}( \adsequence, D) \in \C^b( \TAb)$ is strictly exact.
\end{prop}
\begin{proof}
  1) Suppose that there is an open subgroup $U \subseteq D^2$ whose inverse image in $D^1$ is zero. Then
  $D^1$ is discrete, and the map $D^1 \oplus U \to D^2$ is an open embedding. Proposition \ref{div_disc}
  allows to extend the projection $D^1 \oplus U \twoheadrightarrow D^1$ to a morphism $D^2 \to D^1$ that
  splits the sequence $D^{\bullet}$. Hence i holds in this case; dually, iii holds if there is a compact
  subgroup $K \subseteq C^2$ which surjects onto $C^3$.

  2) The functor $\Hom( \_, D)$ is left exact; in order to prove iii, we have to show that $\iota: C^1
  \hookrightarrow C^2$ induces a strict epimorphism $\iota^*: \Hom( C^2, D) \to \Hom( C^1, D)$. Since $C_{
  \torus}^{\bullet} = 0$, proposition \ref{types_sequence} implies that $C_{\reals}^{\bullet} \subseteq C^{
  \bullet}$ is a split exact subcomplex; thus it is a direct summand in $\C^b(\Codiv)$ because $\reals$ is
  injective in $\LCA$ according to \cite[Thm. 3.2]{moskowitz}. Splitting it off, we may assume $C^2_{\reals}
  = 0$; then we can find a compact open subgroup $V \subseteq C^2$. We factor $\iota: C^1 \hookrightarrow C^2
  $ into $\iota_1: C^1 \hookrightarrow V + \iota( C^1)$ followed by the open inclusion $\iota_2: V + \iota(
  C^1) \hookrightarrow C_2$; then $\iota_1^*$ is a strict epimorphism by step 1, and $\iota_2^*$ is so if $D$
  is \emph{strictly} divisible by proposition \ref{div_disc}.

  3) Suppose that $D$ is discrete or $D = \reals$ or $D = \torus$. Then iii holds by step 2, and iv also
  holds: It obviously holds for $D = \reals$ and $D = \torus$, and if $D$ is discrete, then it is torsion by
  the hypothesis of iv; hence $\Hom( \integers, D) \cong D \cong \Hom( \reals \oplus \hat{\integers}, D)$
  topologically. Since $\rationals \to \adeles$ is a pushout of the natural map $\integers \to \reals \oplus
  \hat{\integers}$ and $\Hom(\_,D)$ is left exact, we get $\Hom(\rationals, D) \cong \Hom( \adeles, D)$
  topologically as well. So iii and iv hold for these $D$; by duality, i and ii follow if $C$ is compact or
  $C = \reals$ or $C = \integers$. \cite[Thm. 2.5]{moskowitz} states that these cover all cases in which $C$
  is \emph{compactly generated}, i.\,e.\ generated as an abstract group by some compact subset.

  4) The general case follows from step 3 by means of lemma \ref{limits}: \cite[Thm. 5.14]{hewittross} states
  that every compact subset of $C$ is contained in a compactly generated open subgroup $\tilde{C} \subseteq C
  $. Thus $\Hom( C, A)$ is the topological projective limit of the $\Hom( \tilde{C}, A)$ for every LCA group
  $A$. Note that $\Hom(C, A)$ surjects onto $\Hom( \tilde{C}, A)$ for divisible $A$ by proposition
  \ref{div_disc}, the kernel being $\Hom( C/\tilde{C}, A_{\disc})$ because $C/\tilde{C}$ is discrete. Since
  $D^1_{\disc}$ and $\rationals$ are injective in $\Ab$, lemma \ref{limits} applies here, completing the
  proof of i and ii. iii and iv follow by duality.
\end{proof}

\begin{cor} \label{acyclic_long}
  If $C^{\bullet} \in \C^b( \Codiv)$ and $D^{\bullet} \in \C^b( \Div)$ satisfy one of the following four
  conditions, then $\Hom^{\bullet}( C^{\bullet}, D^{\bullet}) \in \C^b( \TAb)$ is strictly exact:

  i) $D^{\bullet}$ is strictly exact with $D_{\integers}^{\bullet} = 0$.

  ii) $D^{\bullet} = \adsequence$, and $C_{\torus}^{\bullet} = 0$.

  iii) $C^{\bullet}$ is strictly exact with $C_{\torus}^{\bullet} = 0$.

  iv) $C^{\bullet} = \adsequence$, and $D_{\integers}^{\bullet} = 0$.
\end{cor}
\begin{proof}
  i and ii) Filtering $C^{\bullet}$ by its stupid truncations and using the resulting long exact cohomology
  sequences in $\LH( \TAb)$, we may assume that $C^{\bullet} = C$ is a single object of $\Codiv$. For i,
  we can furthermore decompose $D^{\bullet}$ into short strictly exact sequences; note that each group
  appearing here is a quotient of some $D^n$ and hence still divisible without $\integers$-part.
  Now the previous proposition applies.

  iii and iv) follow dually.
\end{proof}

\begin{prop} \label{resolve}
  Let $A$ be an LCA group.

  i) There is an open embedding $A \hookrightarrow D$ with $D \in \Div$ and $D/A$ torsion.

  ii) If $A \in \FLCA$, then we can achieve $D \in \FLCA$ in i.

  iii) If $A \in \P$, then we can achieve $D = \rationals^r \oplus \reals^s \oplus \textresprod_p (
  \rationals_p^{r_p} : \integers_p^{r_p})$ in i.
\end{prop}
\begin{proof}
  1) There always is a divisible abelian group $D$ containing $A_{\disc}$ such that $D/A_{\disc}$ is torsion,
  cf. \cite[\S24]{fuchs}. We call a subset in $D$ open if and only if its intersection with each coset modulo
  $A$ is open in $A$; this makes $D$ an LCA group containing $A$ as an open subgroup, thereby proving i.

  2) Suppose that $A \in \FLCA$ is a topological $p$-group. Then proposition \ref{p_structure} reduces ii to
  the special cases $A = \integers/p^n$, $\rationals_p/\integers_p$ or $\integers_p$, $\rationals_p$, in
  which $D = \rationals_p/\integers_p$ or $\rationals_p$ does the trick. Under the hypothesis of iii, $A$ is
  torsionfree; hence $\integers_p^r \subseteq A \subseteq \rationals_p^r$ for some $r$ by proposition
  \ref{p_structure}, so $D := \rationals_p^r$ proves iii here.
   
  3) Suppose that $A \in \FLCA$ is of type $\adeles$. Writing $A_{\adeles} \cong A_{\reals} \oplus
  \textresprod_p (A_p: U_p)$, let $A_p \hookrightarrow D_p$ be the open embedding constructed in the previous
  step 2. Then $D := A_{\reals} \oplus \textresprod_p (D_p : U_p)$ proves ii. If $A \in \P$, then this $D$
  also proves iii: It has the required form because any compact open subgroup $U \subseteq \rationals_p^r$
  differs from the standard subgroup $\integers_p^r \subseteq \rationals_p^r$ only by an automorphism of
  $\rationals_p^r$.

  4) Now suppose that $A \in \FLCA$ is arbitrary. Let $\iota': A_{\torus} \to D' := A_{\torus}$ be the
  identity, let $\iota'': A_{\adeles} \hookrightarrow D''$ be the open embedding constructed in step 3, and
  let $\iota''': A_{\integers} \hookrightarrow D''':= A_{\integers} \otimes \rationals$ be the canonical map.
  Since $D'$, $D''$ and $D'''$ are divisible, there is a morphism $\iota: A_{\disc} \to D := (D' \oplus D''
  \oplus D''')_{\disc}$ in $\Ab$ that respects the obvious three step filtrations and induces $\iota_{\disc}'
  $, $\iota_{\disc}''$, $\iota_{\disc}'''$ on the filtration subquotients. Since these are injective and
  their cokernels are torsion groups that have finite ranks, the same holds for $\iota$ due to proposition
  \ref{fin_rank_ext}. We endow $D$ with the unique group topology for which $\iota$ is an open embedding,
  cf. step 1; then $D \in \FLCA$ by proposition \ref{fin_rank_ext}, which completes the proof of ii.

  5) Finally, let $A \in \P$ be arbitrary. Let $\iota': A_{\adeles} \hookrightarrow D'$ be the open embedding
  constructed in step 3, and let $\iota'': A_{\integers} \hookrightarrow D'' := A_{\integers} \otimes
  \rationals$ be the canonical map. Due to proposition \ref{div_disc}, we can extend $\iota'$ to a morphism
  $A \to D'$; its direct sum $\iota: A \to D := D' \oplus D''$ with $\iota'' \circ \pi: A \twoheadrightarrow
  A_{\integers} \hookrightarrow D''$ is injective with torsion cokernel because $\iota'$ and $\iota''$ are,
  and $\iota$ is open because its restriction $\iota'$ to the open subgroup $A_{\adeles} \subseteq A$ is.
  Since $D$ has the required form, this proves iii.
\end{proof}

\begin{cor} \label{fin_resolve}
  Every $A \in \FLCA$ admits a closed embedding $A \hookrightarrow I$ with $I \in \I$.
\end{cor}
\begin{proof}
  Using proposition \ref{resolve}.ii, we may assume that $A$ is divisible. Then corollary \ref{div_split}
  yields $A \cong F_{\integers} A \oplus \rationals^r$ for some $r$, and we can simply take $I := F_{
  \integers} A \oplus \adeles^r$.
\end{proof}

\begin{cor} \label{resolve_long}
  i) Every bounded complex $A^{\bullet} \in \C^b( \LCA)$ admits a strict quasi-isomorphism $f: A^{\bullet}
  \to D^{\bullet}$ with $D^{\bullet} \in \C^b( \Div)$.

  ii) The class $\Quis$ of strict quasi-isomorphisms in $\K^b( \Div)$ is localizing. The resulting category
  $\D^b(\Div) := \K^b(\Div)/\Quis$ is triangulated, and the inclusion functor $\Div \hookrightarrow \LCA$
  induces a triangulated equivalence $\D^b( \Div) \to \D^b( \LCA)$.

  iii) i and ii remain true if $\Div \subseteq \LCA$ is replaced by $\I \subseteq \FLCA$.
\end{cor}
\begin{proof}
  i) We construct $f^n: A^n \to D^n$ inductively; as $A^n = 0$ for $n \ll 0$, we can start with $D^n = 0$
  for $n \ll 0$. Suppose that $\ldots D^{n-1} \to D^n$ and $\ldots f^{n-1}, f^n$ are already constructed such
  that the mapping cone of $f$ has a strict boundary operator $F^n: A^n \oplus D^{n-1} \to A^{n+1} \oplus D^n
  $ with $\ker( F^n) = \im( F^{n-1})$. Using proposition \ref{resolve}.i, we can find a strict monomorphism
  $\coker( F^n) \to D^{n+1}$ with $D^{n+1} \in \Div$; the components of the composition $A^{n+1} \oplus D^n
  \twoheadrightarrow \coker( F^n) \hookrightarrow D^{n+1}$ yield the required morphisms $f^{n+1}: A^{n+1} \to
  D^{n+1}$ and $\partial: D^n \to D^{n+1}$. This constructs a strict quasi-isomorphism $f: A^{\bullet} \to 
  D^{\bullet}$. If $A^{n+1} = A^{n+2} = 0$, then $\coker( F^n)$ is already divisible, so we can take $D^{n+1}
  = \coker( F^n)$ and $D^{n+2} = 0$; thus we can arrange that $D^{\bullet}$ is bounded above and hence in
  $\C^b( \Div)$.

  ii) Because the class of strictly exact complexes is stable under extensions in $\K^b( \LCA)$, it is so in
  the triangulated category $\K^b( \Div)$ as well; this implies that $\Quis$ is localizing in $K^b( \Div)$
  and compatible with the triangulation, so the localized category $\D^b( \Div)$ inherits a triangulation by
  \cite[Thm. IV.2.2]{manin}. The inclusion functor $\Div \hookrightarrow \LCA$ induces a functor of
  triangulated categories $\D^b( \Div) \to \D^b( \LCA)$ which is essentially surjective by i and fully
  faithful by \cite[Prop. III.2.10]{manin} (whose hypothesis $b_2$ holds here by i again).

  iii) same proof, using corollary \ref{fin_resolve} instead of proposition \ref{resolve}.i.
\end{proof}

\begin{rem}
  These results can be dualized in a straightforward way. For example, the class of strict quasi-isomorphisms
  is also localizing in $\K^b( \Codiv)$ and in $\K^b( \P)$, leading to derived categories $\D^b( \Codiv)$ and
  $\D^b( \P)$ which are equivalent to $\D^b( \LCA)$ and to $\D^b( \FLCA)$, respectively.
\end{rem}

\begin{prop}
  If $A, B \in \FLCA$, then $\Hom( A, B) \in \FLCA$.
\end{prop}
\begin{proof}
  Using corollary \ref{fin_resolve}, left exactness of $\Hom( A, \_)$ and proposition \ref{fin_rank_ext}, we
  may assume $B \in \I$ without loss of generality; dually, we may furthermore assume $A \in \P$. Now $B$ is
  divisible and has finite ranks, so $n \cdot \_: B \to B$ is a strict epimorphism whenever $n$ is prime;
  hence it is so whenever $n$ is a product of primes, i.\,e.\ $B$ is strictly divisible. Using propositions
  \ref{resolve}.iii, \ref{div_disc} and \ref{fin_rank_ext}, we may thus assume $A = \rationals$, $\reals$ or
  $\textresprod_p (\rationals_p^{r_p}: \integers_p^{r_p})$; dually, we may furthermore assume $B = \solenoid
  $, $\reals$ or $\textresprod_p(\rationals_p^{s_p}: \integers_p^{s_p})$. These special cases can be
  checked explicitly.
\end{proof}

\begin{defn} \label{tensor}
  The \emph{tensor product} of LCA groups $A, B \in \FLCA$ is the LCA group
  $A \otimes B := \Hom( A, B^{\dual})^{\dual} \in \FLCA$.
\end{defn}

\begin{prop} \label{monoidal}
  $(\FLCA, \otimes, \Hom)$ is a closed symmetric monoidal category.
\end{prop}
\begin{proof}
  \cite[Thm. 4.2]{moskowitz} implies that $A_1 \otimes (A_2 \otimes ( \ldots (A_{n-1} \otimes A_n) \ldots ))$
  is for any $A_1, \ldots, A_n \in \FLCA$ canonically isomorphic to the dual of the group of all continuous
  multiadditive maps $A_1 \times \ldots \times A_n \to \torus$, endowed with the compact-open topology.
  Thus we see that $(\FLCA, \otimes)$ is a symmetric monoidal category with unit object $\integers$; cf.
  \cite[Chapter VII]{maclane} for the categorical terminology. \cite[Thm. 4.2]{moskowitz} also implies that
  $\Hom( A \otimes B, C)$ is canonically isomorphic to $\Hom( A, \Hom( B, C))$ for $A, B, C \in \FLCA$, so
  our monoidal category is closed.
\end{proof}

\section{Derived $\Hom$-functors} \label{derived_Hom}

\begin{prop} \label{RHom_fin}
  The left exact bifunctor $\Hom: \FLCA^{\op} \times \FLCA \to \FLCA$ has a right derived functor $\RHom:
  \D^b( \FLCA)^{\op} \times \D^b( \FLCA) \to \D^b( \FLCA)$.
\end{prop}
\begin{proof}
  This follows from \ref{acyclic} and \ref{fin_resolve} by standard methods; cf. \cite[Thm. III.6.8]{manin}.
  More precisely, corollary \ref{acyclic_long} implies that $\Hom^{\bullet}: \K^b( \P)^{\op} \times \K^b( \I)
  \to \K^b( \FLCA)$ induces a bifunctor of triangulated categories $\D^b( \P)^{\op} \times \D^b( \I) \to
  \D^b( \FLCA)$. Given $A^{\bullet}, B^{\bullet} \in \C^b( \FLCA)$, we use corollary \ref{resolve_long} and
  its dual to choose strict quasi-isomorphisms $r: P^{ \bullet} \to A^{ \bullet}$ and $c: B^{\bullet} \to
  I^{\bullet}$ with $P^{\bullet} \in \C^b( \P)$ and $I^{\bullet} \in \C^b( \I)$; then we define $\RHom(
  A^{\bullet}, B^{\bullet}) := \Hom^{\bullet}( P^{\bullet}, I^{\bullet})$. This is a well-defined bifunctor
  of triangulated categories $\D^b( \FLCA)^{\op} \times \D^b( \FLCA) \to \D^b( \FLCA)$ by corollary
  \ref{resolve_long} and its dual; the functorial morphism $r^* c_*: \Hom^{\bullet}(A^{\bullet}, B^{\bullet})
  \to \RHom(A^{\bullet}, B^{\bullet})$ clearly has the required universal property by construction.
\end{proof}

  \begin{sidewaystable}
  \begin{equation*} \renewcommand{\arraystretch}{1.5}\begin{array}{|c||c|ccc|ccc|ccc|}
    \hline $\backslashbox{A}{B}$            & \integers/p^n
      & \rationals_p/\integers_p     & \rationals_p     & \integers_p
      & \torus                       & \reals           & \integers
      & \solenoid                    & \adeles          & \rationals\\\hline \hline
    \integers/p^n                    & \integers/p^n \to[0] \integers/p^n
      & \integers/p^n                & 0                & \integers/p^n[-1]
      & \integers/p^n                & 0                & \integers/p^n[-1]
      & 0                            & 0                & 0\\\hline
    \integers_p                      & \integers/p^n
      & \rationals_p/\integers_p     & \rationals_p     & \integers_p
      & \rationals_p/\integers_p     & 0                & \rationals_p/\integers_p[-1]
      & \rationals_p                 & \rationals_p     & 0\\
    \rationals_p                     & 0
      & \rationals_p                 & \rationals_p     & 0
      & \rationals_p                 & 0                & \rationals_p[-1]
      & \rationals_p                 & \rationals_p     & 0\\
    \rationals_p/\integers_p         & \integers/p^n[-1]
      & \integers_p                  & 0                & \integers_p[-1]
      & \integers_p                  & 0                & \integers_p[-1]
      & 0                            & 0                & 0\\\hline
    \integers                        & \integers/p^n
      & \rationals_p/\integers_p     & \rationals_p     & \integers_p
      & \torus                       & \reals           & \integers
      & \solenoid                    & \adeles          & \rationals\\
    \reals                           & 0
      & 0                            & 0                & 0
      & \reals                       & \reals           & 0
      & \reals                       & \reals           & 0\\
    \torus                           & \integers/p^n[-1]
      & \rationals_p/\integers_p[-1] & \rationals_p[-1] & \integers_p[-1]
      & \integers                    & 0                & \integers[-1]
      & \rationals \into \finadeles  & \finadeles[-1]   & \rationals[-1]\\\hline
    \rationals                       & 0
      & \rationals_p                 & \rationals_p     & 0
      & \solenoid                    & \reals           & \rationals \into \finadeles 
      & \solenoid                    & \adeles          & \rationals\\
    \adeles                          & 0
      & \rationals_p                 & \rationals_p     & 0
      & \adeles                      & \reals           & \finadeles[-1]
      & \adeles                      & \adeles          & 0\\
    \solenoid                        & 0
      & 0                            & 0                & 0
      & \rationals                   & 0                & \rationals[-1]
      & \rationals                   & 0                & \rationals[-1]\\ \hline
  \end{array} \end{equation*}
    \caption{The complex $\RHom(A,B) \in \D^b(\FLCA)$ for various LCA groups $A, B \in \FLCA$} \label{table}
  \end{sidewaystable}
\begin{examples} \label{RHom_examples}
  Table \ref{table} lists $\RHom( A, B)$ for various $A, B \in \FLCA$. Furthermore, it is easy to check
  that $\RHom( A, B) = 0$ if $A$ is a topological $p$-group and $B$ is a topological $q$-group for prime
  numbers $p \neq q$. 
\end{examples}

\begin{rem}
  i) Dually, the right exact bifunctor $\otimes: \FLCA \times \FLCA \to \FLCA$ of \ref{tensor} has a left
  derived functor $\Lotimes: \D^b( \FLCA) \times \D^b( \FLCA) \to \D^b( \FLCA)$; it is given by
  $A^{\bullet} \Lotimes B^{\bullet} := \RHom( A^{\bullet}, (B^{\bullet})^{\dual})^{\dual}$.

  ii) Proposition \ref{monoidal} implies that $(\D^b( \FLCA), \Lotimes, \RHom)$ is a closed
  symmetric monoidal category as well.
\end{rem}

We denote by $\mathrm{K}_0( \FLCA)$ the abelian group generated by symbols $[A]$ for each $A \in \FLCA$
subject to the relations $[B] = [A] + [C]$ for all strictly exact sequences $0 \to A \to B \to C \to 0$;
this abelian group comes with an automorphism of order $2$, given by $[A] \mapsto [A^{\dual}]$. We write
$[A^{\bullet}] := \sum_n (-1)^n [A^n] \in \mathrm{K}_0( \FLCA)$ for complexes $A^{\bullet} \in \C^b(\FLCA)$.
Then $\mathrm{K}_0( \FLCA)$ becomes a commutative ring if we put $[A^{\bullet}] \cdot [B^{ \bullet}] :=
[A^{\bullet} \Lotimes B^{\bullet}]$; this is easily checked to be well defined. The ring structure is not
compatible with the duality involution; explicitly, we have:

\begin{prop}
  Let $v$ run over all places of $\rationals$, i.\,e.\ $v = p$ is a prime or $v = \infty$.

  i) Sending $(r_v, s_v)_v \in \prod_v \naturals^2$ to $r_{\infty}[ \integers] + s_{\infty}[ \torus] +
  [\prod_p \integers_p^{r_p}] + [\bigoplus_p (\rationals_p/\integers_p)^{s_p}]$ defines a group isomorphism
  $\prod_v \integers^2 \to \mathrm{K}_0( \FLCA)$ under which the involution $[A] \mapsto [A^{\dual}]$ on
  $\mathrm{K}_0(\FLCA)$ corresponds to the involution $(r_v, s_v)_v \mapsto (s_v, r_v)_v$ on $\prod_v
  \integers^2$.

  ii) Sending $(r_v, s_v)_v \in \prod \naturals^2$ to $r_{\infty}[ \reals] - s_{\infty}[ \solenoid] +
  [\textresprod_p (\rationals_p^{r_p}: \integers_p^{r_p})] - [\bigoplus_p (\rationals_p/ \integers_p)^{s_p}]$
  instead defines a ring isomorphism $\prod_v \integers^2 \to \mathrm{K}_0( \FLCA)$.
\end{prop}
\begin{proof}
  i) The given map on $\prod_v \naturals^2$ extends canonically to a group homomorphism on $\prod_v \integers
  ^2$ which is obviously compatible with the involutions in question. Its image contains all divisible
  discrete torsion groups by lemma \ref{p_structure}, so it contains all objects of $\P$ due to proposition
  \ref{resolve}.iii; hence this map is surjective according to the dual of corollary \ref{fin_resolve}.
  In order to prove injectivity, we construct a left inverse by sending the class $[A] \in K_0( \FLCA)$ of
  $A \in \FLCA$ to the integers $r_{\infty} := \dim_{\reals} \Hom( A, \reals)$ and $s_{\infty} := \dim_{
  \reals} \Hom( \reals, A)$ and $r_p := s_{\infty} + \sum_n (-1)^n \dim_{\rationals_p} \H^n(\RHom(A,
  \rationals_p))$ and $s_p:= r_{\infty} + \sum_n (-1)^n \dim_{\rationals_p} \H^n( \RHom( \rationals_p, A))$;
  that these are finite-dimensional vector spaces can be checked on generators of $K_0( \FLCA)$, and it is
  easily verified for the generators that we have just obtained by proving surjectivity.

  ii) Using the examples \ref{RHom_examples}, it is easy to see that this determines a well-defined ring
  homomorphism; its bijectivity can be deduced from i.
\end{proof}

Our next aim is to extend the derived functor $\RHom$ to all of $\D^b( \LCA)$. Here the main problem is the
lack of enough acyclic objects, cf. \cite[Thm. 3.6]{moskowitz}; instead, we will use divisible and
codivisible groups -- which are `almost acyclic' by proposition \ref{acyclic} -- and the `standard
resolution' $0 \to \rationals \to \adeles \to \solenoid \to 0$.
\begin{defn} \label{RHom_def}
  Given complexes $C^{\bullet} \in \C^b( \Codiv)$ and $D^{\bullet} \in \C^b( \Div)$, we define $\RHom( C^{
  \bullet}, D^{\bullet}) \in \C^b( \TAb)$ to be the mapping cone of the composition pairing
  \begin{equation*}
    \Hom^{\bullet}(\adsequence, D^{\bullet})_{\disc} \otimes
    \Hom^{\bullet}(C^{\bullet}, \adsequence)_{\disc} \longto[\circ] \Hom^{\bullet}(C^{\bullet}, D^{\bullet})
  \end{equation*}
  where the complex $\adsequence \in \C^b( \LCA)$ is located in degrees $-1$, $0$ and $1$.
\end{defn}
Here the tensor product is just one of complexes of discrete rational vector spaces (because $\rationals$,
$\adeles$ and $\solenoid$ are topological rational vector spaces).

Clearly, this $\RHom$ is an additive bifunctor $\C^b( \Codiv) \times \C^b( \Div) \to \C^b( \TAb)$ which
respects homotopies, shifts and mapping cones; thus it descends to a bifunctor of triangulated categories
$\K^b( \Codiv) \times \K^b( \Div) \to \K^b( \TAb)$. Furthermore, one has a natural duality isomorphism
$\RHom( D^{\bullet \dual}, C^{\bullet \dual}) \cong \RHom( C^{\bullet}, D^{\bullet})$.
\begin{lemma} \label{disc_equiv}
  The canonical restriction and projection morphisms
  \begin{align*}
    \rho: \Hom^{\bullet}( C^{\bullet}, \adsequence) & \longto
      \Hom^{\bullet}( C_{\torus}^{\bullet}, \adsequence) \cong (C_{\torus}^{\bullet})^{\dual}[-1],\\
    \pi: \Hom^{\bullet}( \adsequence, D^{\bullet}) & \longto
      \Hom^{\bullet}( \adsequence, D_{\integers}^{\bullet}) \cong D_{\integers}^{\bullet}[-1]
  \end{align*}
  are strict quasi-isomorphisms for all $C^{\bullet} \in \C^b( \Codiv)$ and $D^{\bullet} \in \C^b( \Div)$.
\end{lemma}
\begin{proof}
  Corollary \ref{div_split}.i implies that $\pi$ is objectwise a split epimorphism, so its mapping cone is
  homotopy equivalent to $\ker( \pi) = \Hom^{\bullet}(\adsequence, F_{\integers} D^{\bullet})$ which is
  strictly exact by corollary \ref{acyclic_long}.iv. The claim about $\rho$ follows dually.
\end{proof}
\begin{cor} \label{triangle}
  One has a canonical distinguished triangle
  \begin{equation*}
    \big( D^{\bullet}_{\integers} \otimes (C^{\bullet}_{\torus})^{\dual} \big)[-2] \longto
    \Hom^{\bullet}( C^{\bullet}, D^{\bullet}) \longto \RHom( C^{\bullet}, D^{\bullet}) \longto
    \big( D^{\bullet}_{\integers} \otimes (C^{\bullet}_{\torus})^{\dual} \big)[-1]
  \end{equation*}
  in $\K^b(\TAb)$ which is functorial in $C^{\bullet} \in \K^b(\Codiv)$ and $D^{\bullet} \in \K^b(\Div)$.
  In particular, $\RHom( C^{\bullet}, D^{\bullet}) \cong \Hom^{\bullet}( C^{\bullet}, D^{\bullet})$ in
  $\K^b(\TAb)$ if $C_{\torus}^{\bullet} = 0$ or $D_{\integers}^{\bullet} = 0$.
\end{cor}
\begin{proof}
  $\pi \otimes \rho$ is a quasi-isomorphism from the tensor product complex in definition \ref{RHom_def} to
  $\big( D^{\bullet}_{\integers} \otimes (C^{\bullet}_{\torus})^{\dual} \big)[-2]$, so it is an isomorphism
  in $\K^b( \TAb)$ by linear algebra.
\end{proof}

\begin{thm} \label{RHom}
  The bifunctor $\RHom$ of definition \ref{RHom_def} induces a bifunctor
  \begin{equation*}
    \RHom: \D^b( \LCA)^{\op} \times \D^b( \LCA) \longto \D^b( \TAb)
  \end{equation*}
   of triangulated categories and a morphism $s: \Hom^{\bullet}( A^{\bullet}, B^{\bullet}) \to \RHom( A^{
   \bullet}, B^{\bullet})$ in $\D^b(\TAb)$ that is functorial in $A^{\bullet}, B^{\bullet} \in \C^b(\LCA)$.
\end{thm}
\begin{proof}
  1) Let $V$ be a discrete rational vector space, say in degree $0$. By its very definition, $\RHom(
  \adsequence, V)$ is the mapping cone of the pairing
  \begin{equation*}
    \circ: V[-1] \otimes \Hom^{\bullet}( \adsequence, \adsequence)_{\disc} \longto V[-1]
  \end{equation*}
  which is a homotopy equivalence because the second factor is homotopy equivalent to $\rationals$ by lemma
  \ref{disc_equiv}; hence $\RHom( \adsequence, V)$ is strictly exact.

  2) Let $D^{\bullet} \in \C^b( \Div)$. Corollary \ref{div_split}.i provides us with a distinguished triangle
  $F_{\integers} D^{\bullet} \to D^{\bullet} \to D_{\integers}^{\bullet} \to F_{\integers} D^{\bullet}[1]$ in
  $\K^b( \Div)$. $\RHom( \adsequence, F_{\integers} D^{\bullet})$ is strictly exact due to the corollaries
  \ref{triangle} and \ref{acyclic_long}.iv; since $D_{\integers}^{\bullet}$ is homotopy equivalent to the
  direct sum of its cohomology by linear algebra, step 1 above implies that $\RHom( \adsequence, D_{\integers
  }^{\bullet})$ is also strictly exact. Using \cite[Prop. 1.2.14]{schneiders}, this shows that
  $\RHom( \adsequence, D^{\bullet})$ is strictly exact as well.

  3) Let $D^{\bullet} \in \C^b(\Div)$ be strictly exact; we claim that for all $C^{\bullet} \in \C^b( \Codiv)
  $, the complex $\RHom(C^{\bullet}, D^{\bullet}) \in \C^b(\TAb)$ is strictly exact. We prove this by
  induction on $\sum_n \dim_{\rationals} \H^n( D_{\integers}^{\bullet})$ which is finite by remark
  \ref{finiteness}.i.

  If this sum is zero, then $F_{\integers} D^{\bullet} \hookrightarrow D^{\bullet}$ is a homotopy equivalence
  by corollary \ref{div_split}.i, so $\RHom( C^{\bullet}, D^{\bullet})$ is strictly exact due to the
  corollaries \ref{triangle} and \ref{acyclic_long}.i.

  For the induction step, suppose that there is a nonzero class $\gamma \in \H^n( D_{\integers}^{\bullet})$.
  Lemma \ref{disc_equiv} implies that there is a morphism $\tilde{\gamma}: \adsequence[-n-1] \to D^{\bullet}$
  in $\C^b( \Div)$ such that the induced morphism $\tilde{\gamma}_{\integers}: \rationals[-n] \to D_{
  \integers}^{\bullet}$ maps $1 \in \rationals$ to a cycle representing $\gamma$. Denoting by ${D'}^{\bullet}
  $ the mapping cone of $\tilde{\gamma}$, we get a distinguished triangle $\rationals[-n] \to[\tilde{\gamma
  }_{\integers}] D_{\integers}^{\bullet} \to D_{\integers}'^{\bullet} \to \rationals[-n+1]$ which shows that
  the induction hypothesis applies to ${D'}^{\bullet}$. Since $\RHom( C^{\bullet}, \adsequence)$ is also
  strictly exact by the dual of step 2 above, \cite[Prop. 1.2.14]{schneiders} completes the induction step.

  4) The previous step 3 and its dual imply that the functor $\RHom$ of \ref{RHom_def} induces a bifunctor
  of triangulated categories $\D^b( \Codiv)^{\op} \times \D^b( \Div) \to \D^b( \TAb)$. Given $A^{\bullet},
  B^{\bullet} \in \C^b( \LCA)$, we use corollary \ref{resolve_long}.i and its dual to choose strict
  quasi-isomorphisms $r: C^{\bullet} \to A^{\bullet}$ and $c: B^{\bullet} \to D^{\bullet}$ with $C^{\bullet}
  \in \C^b( \Codiv)$ and $D^{\bullet} \in \C^b( \Div)$; then we define $\RHom( A^{\bullet}, B^{\bullet}) :=
  \RHom^{\bullet}( C^{\bullet}, D^{\bullet})$. This is a well-defined bifunctor of triangulated categories
  $\D^b( \LCA)^{\op} \times \D^b( \LCA) \to \D^b( \TAb)$ by corollary \ref{resolve_long}.ii and its dual;
  the required functorial morphism is given by $s := r^* c_*: \Hom^{\bullet}( A^{\bullet}, B^{\bullet}) \to
  \RHom( A^{\bullet}, B^{\bullet})$.
\end{proof}

\begin{rem}
  The bifunctor $\RHom$ of theorem \ref{RHom} actually extends the bifunctor $\RHom$ of proposition
  \ref{RHom_fin}; more precisely, both induce the same bifunctor $\D^b( \FLCA)^{\op} \times \D^b( \FLCA) \to
  \D^b(\TAb)$. This follows from the fact that both induce the same bifunctor $\K^b( \P)^{\op} \times \K^b(
  \I) \to \D^b( \TAb)$ due to corollary \ref{triangle}.
\end{rem}

\begin{example} \label{RHom_example}
  For arbitrary index sets $I$ and $J$,
  \begin{equation*}
    \RHom( \prod_I \torus, \bigoplus_J \integers) \cong \RHom( \prod_I \hat{\integers} \to \prod_I \solenoid,
    \bigoplus_J \rationals \to \bigoplus_J \rationals/\integers) \cong \bigoplus_{I \times J} \integers[-1]
  \end{equation*}
  where the sums are discrete and the products carry the Tychonoff topology.
\end{example}

\begin{defn}
  $\Ext^n( A^{\bullet}, B^{\bullet}) := \H^n( \RHom( A^{\bullet}, B^{\bullet})) \in \LH( \TAb)$ for $n \in
  \integers$ and bounded complexes $A^{\bullet}$ and $B^{\bullet}$ of LCA groups.
\end{defn}
Note that $\Ext^n: \D^b( \LCA)^{\op} \times \D^b( \LCA) \to \LH( \TAb)$ and the composed functor $\Ext^n_{
\disc}: \D^b( \LCA)^{\op} \times \D^b( \LCA) \to \LH( \TAb) \to \Ab$ are cohomological bifunctors, i.\,e.
fixing one variable, they transform distinguished triangles in the other variable to long exact sequences.
Moreover, we have a canonical duality isomorphism $\Ext^n(B^{\bullet \dual}, A^{\bullet \dual}) \cong \Ext^n(
A^{\bullet}, B^{\bullet})$ inherited from $\RHom$.
\begin{prop} \label{RHom_lim}
  For $A^{\bullet}, B^{\bullet} \in \C^b( \LCA)$, there is a canonical isomorphism

  i) $\Ext^n_{\disc}(A^{\bullet}, B^{\bullet}) \cong \Hom_{\D^b(\LCA)}(A^{\bullet}, B^{\bullet}[n])$ in $\Ab$

  ii) $\Ext^n( A^{\bullet}, B^{\bullet}) \cong \varinjlim \H^n( \Hom^{\bullet}( {A'}^{\bullet}, {B'}^{\bullet
  }))$ in $\LH( \TAb)$, the limit being over all strict quasi-isomorphisms ${A'}^{\bullet} \to A^{\bullet}$
  and $B^{\bullet} \to {B'}^{\bullet}$ in $\C^b( \LCA)$.
\end{prop}
\begin{proof}
  Due to corollary \ref{resolve_long}.i and its dual, we may assume $A^{\bullet} \in \C^b( \Codiv)$ and
  $B^{\bullet} \in \C^b( \Div)$ without loss of generality, and it suffices to consider strict
  quasi-isomorphisms with ${A'}^{\bullet} \in \C^b( \Codiv)$ and ${B'}^{\bullet} \in \C^b( \Div)$ in ii.

  Given a cohomology class $\gamma \in \H^p( B_{\integers}'^{\bullet})$, lemma \ref{disc_equiv} implies that
  there is a morphism $\tilde{\gamma}: \adsequence[-p-1] \to {B'}^{\bullet}$ in $\C^b( \Div)$ such that the
  induced morphism $\tilde{\gamma}_{\integers}: \rationals[-p] \to B_{\integers}'^{\bullet}$ maps $1 \in
  \rationals$ to a cycle representing $\gamma$. Let ${B''}^{\bullet}$ be the mapping cone of $\tilde{\gamma}
  $; then the natural strict quasi-isomorphism ${B'}^{\bullet} \to {B''}^{\bullet}$ maps $\gamma$ to $0 \in
  \H^p( B_{\integers}''^{\bullet})$ by construction. This shows that the inductive limit of the $\H^p( B_{
  \integers}'^{\bullet})$ in $\Ab$ vanishes; hence the inductive limit of the $\H^n( B_{\integers}'^{\bullet}
  \otimes (A_{\torus}'^{\bullet})^{\dual})$ in $\Ab$ also vanishes. By corollary \ref{triangle} and the
  exactness of inductive limits in $\Ab$,
  \begin{equation*}
    \varinjlim s_{\disc}: \varinjlim \H^n( \Hom^{\bullet}( {A'}^{\bullet}, {B'}^{\bullet}))_{\disc} \longto
    \varinjlim \Ext^n({A'}^{\bullet}, {B'}^{\bullet})_{\disc}
  \end{equation*}
  is thus an isomorphism of abelian groups. Here $\H^n( \Hom^{\bullet}( {A'}^{\bullet}, {B'}^{\bullet}))_{
  \disc}$ is precisely the group of morphisms ${A'}^{\bullet} \to {B'}^{\bullet}[n]$ in $\K^b( \LCA)$, so
  their inductive limit is the group of morphisms $A^{\bullet} \to B^{\bullet}[n]$ in $\D^b( \LCA)$; this
  implies i. Furthermore,
  \begin{equation*} \xymatrix{
    \H^n( \Hom^{\bullet}( {A'}^{\bullet}, {B'}^{\bullet}))_{\disc} \ar[r]^-{s_{\disc}} \ar[d]
      & \Ext^n({A'}^{\bullet}, {B'}^{\bullet})_{\disc} \ar[d]\\
    \H^n( \Hom^{\bullet}( {A'}^{\bullet}, {B'}^{\bullet}))         \ar[r]^-s
      & \Ext^n({A'}^{\bullet}, {B'}^{\bullet})
  } \end{equation*}
  is a (pullback and) pushout square in $\LH( \TAb)$ because $s$ and $s_{\disc}$ have the same kernel and
  cokernel in $\LH( \TAb)$: This follows from the fact that kernel and cokernel of $s$ in $\LH( \TAb)$
  are discrete groups by corollary \ref{triangle}. Using this pushout property, the inductive limit property
  in question carries over from $\Ab$ to $\LH( \TAb)$.
\end{proof}

\begin{cor}
  $\Ext^n: \D^b( \LCA)^{\op} \times \D^b( \LCA) \to \LH( \TAb)$ is a right derived cohomological functor
  for $\H^n \circ \Hom^{\bullet}: \K^b( \LCA)^{\op} \times \K^b( \LCA) \to \LH( \TAb)$ in the sense of
  \cite[Ch. II, \S 2, D\'{e}f. 1.4]{verdier}.
\end{cor}

In particular, we can consider LCA groups $A$ and $B$ as complexes concentrated in degree zero, obtaining
objects $\Ext^n( A, B)$ of $\LH( \TAb)$. They vanish by construction for $n < 0$, and $\Ext^0( A, B)$ is
canonically isomorphic to $\Hom( A, B)$ because the $\Hom$-functor is left exact. For $n \geq 1$, the
abelian groups $\Ext^n(A, B)_{\disc}$ coincide with the Yoneda-$\Ext^n$-group studied in \cite{fulpgriffith};
this follows from proposition \ref{RHom_lim}.i by standard arguments, cf. \cite[XI.4]{iversen}. Part iv of
the following vanishing result refines \cite[part II, Thm. 2.9]{fulpgriffith}; cf. also
\cite[Section 6]{armacost} for related results.
\begin{prop} \label{vanishing}
  Let $A$, $B$ be LCA groups and $n \geq 1$. Then $\Ext^n( A, B) = 0$ in $\LH( \TAb)$ in each of the
  following cases: 

  i) $A = \reals^n$ or $A = \bigoplus_{j \in J} \integers$ for some index set $J$.

  ii) $A$ is compact codivisible, and $B_{\integers} = 0$.

  iii) $A$ is codivisible with $A_{\torus} = 0$, and $B$ is divisible.

  iv) $n \geq 2$.

  v) $A$ is codivisible, and $B$ is divisible with $B_{\integers} = 0$.

  vi) $A_{\torus} = 0$, and $B$ is discrete divisible.

  vii) $B = \reals^n$ or $B = \prod_J \torus$ for some index set $J$.
\end{prop}
\begin{proof}
  iii and v are consequences of corollary \ref{triangle}.

  Proposition \ref{resolve}.i yields a strictly exact sequence $0 \to B \to D \to D' \to 0$ with $D, D' \in
  \Div$ and $D'$ a discrete torsion group; dually, there also is a strictly exact sequence $0 \to C' \to C
  \to A \to 0$ with $C, C' \in \Codiv$ and $C'$ profinite.

  i) Here $\RHom( A, B)$ is given by the complex $\Hom( A, D) \to \Hom( A, D')$. If $A = \reals^n$, then
  $\Hom( A, D') = 0$; if $A = \bigoplus_{j \in J} \integers$, then $\Hom( A, D) \to \Hom( A, D')$ is
  the Tychonoff product indexed by $J$ of copies of the open surjection $D \to D'$ and thus also an open
  surjection. This proves i and by duality also vii.

  vi) Suppose first that $A$ is a topological torsion group; then $C$ and $C'$ also are, so $\Hom( C,
  \rationals/\integers) \to \Hom( C', \rationals/\integers)$ is surjective by Pontryagin duality. Since $C'$
  is compact and $B$ is discrete, $\Hom( C', B)$ is discrete, and every morphism $C' \to B$ has finite image;
  it thus factors through finitely many summands of $B \cong \bigoplus_j \rationals_{p_j}/\integers_{p_j}$
  \cite[Thm. A.14]{hewittross}, i.\,e.\ through some morphism $(\rationals/\integers)^n \to B$. This shows
  that $\Hom( C, B) \to \Hom( C', B)$ is a strict epimorphism; since this complex computes $\RHom( A, B)$,
  vi follows here.

  Now let $A$ be arbitrary with $A_{\torus} = 0$. Recall that $A_{\reals}$ is a direct summand of $A$;
  choosing a subset of $A$ whose image in the rational vector space $A_{\integers} \otimes \rationals$ is a
  basis, we can construct a strictly exact sequence $0 \to A' \to A \to A'' \to 0$ such that $A' \cong
  \reals^n \oplus \bigoplus_{j \in J} \integers$ and $A''$ is a topological torsion group. We have just seen
  that $\Ext^n( A'', B) = 0$, and $\Ext^n( A', B) = 0$ by i; thus the long exact $\Hom$-$\Ext$-sequence in
  $\LH( \TAb)$ completes the proof of vi. ii follows by duality.

  iv) If $B_{\integers} = 0$, then $D_{\integers} = 0$ as well, so $\RHom( C, B)$ is given by the complex
  $\Hom( C, D) \to \Hom( C, D')$; thus $\Ext^n( C, B) = 0$ for $n \geq 2$. As $\Ext^{n-1}( C', B)$ also
  vanishes by ii, the long exact sequence in $\LH(\TAb)$ implies iv in this case.

  For general $B$, we construct a strictly exact sequence $0 \to B' \to B \to B'' \to 0$ with $B_{\integers
  }'' = 0$ and $B' \cong \bigoplus \integers$ (by choosing a subset of $B$ whose image in the rational vector
  space $B_{\integers} \otimes \rationals$ is a basis); thus the long exact sequence in $\LH(\TAb)$ reduces
  us to the case $B \cong \bigoplus_{j \in J} \integers$. The dual argument allows us to assume $A \cong
  \prod_{i \in I} \torus$ as well; in this case, iv follows from the explicit example \ref{RHom_example}.
\end{proof}

\begin{rem}
  The functors $\Ext^n: \LCA^{\op} \times \LCA \to \LH( \TAb)$ extend canonically from $\LCA$ to $\LH(
  \LCA)$ because the embedding $\LCA \hookrightarrow \D^b( \LCA)$ does so, by the very construction of $\LH(
  \LCA)$. However, part iv of the previous proposition does not extend to $\LH( \LCA)$; its cohomological
  dimension is not $1$, but $2$. For example, let $J$ be an infinite set, and let $A \in \LH( \LCA)$ be the
  cokernel of the natural map $\bigoplus_J \integers/p \hookrightarrow \prod_J \integers/p$, using the
  discrete topology for the direct sum and the Tychonoff topology for the product. Then $\Ext^2( A,
  \integers/p) \neq 0$, e.\,g. because $\Ext^1( \prod_J \integers/p, \integers/p) \cong \bigoplus_J
  \integers/p$ and $\Ext^1( \bigoplus_J \integers/p, \integers/p) \cong \prod_J \integers/p$.
\end{rem}

\end{document}